\documentclass[10pt]{amsart}
\usepackage{amsthm, amsfonts, amssymb, amsmath, amscd}
\usepackage{graphicx}
\usepackage{lmodern}
\usepackage[margin=3.5cm]{geometry}

\newcommand{\To}{\rightarrow}

\newcommand{\into}{\hookrightarrow}
\newcommand{\onto}{\twoheadrightarrow}
\newcommand{\noin}{\noindent}

\newcommand{\C}{\mathbb{C}}
\newcommand{\R}{\mathbb{R}}

\newcommand{\Z}{\mathbb{Z}}
\newcommand{\N}{\mathbb{N}}
\newcommand{\Cstar}{\mathrm{C^*}}
\newcommand{\Cred}{\mathrm{C^*_r}}
\newcommand{\G}{\Gamma}

\newcommand{\M}{\mathrm{M}}
\newcommand{\GL}{\mathrm{GL}}
\newcommand{\Idem}{\mathrm{Idem}}
\newcommand{\tsr}{\mathrm{tsr}}
\newcommand{\hsr}{\mathrm{hsr}}
\newcommand{\bsr}{\mathrm{bsr}}
\newcommand{\csr}{\mathrm{csr}}
\newcommand{\diag}{\mathrm{diag}}
\newcommand{\spec}{\mathrm{sp}}
\newcommand{\hol}{\mathcal{O}}
\newcommand{\supp}{\mathrm{supp}}

\newcommand{\id}{\mathrm{id}}
\newcommand{\Lg}{\mathrm{Lg}}
\newcommand{\re}{\mathrm{Re}}

\newtheorem{thm}{Theorem}[section]
\newtheorem{lem}[thm]{Lemma}
\newtheorem{cor}[thm]{Corollary}
\newtheorem{prop}[thm]{Proposition}

\theoremstyle{definition}
\newtheorem{defn}[thm]{Definition}
\newtheorem{rem}[thm]{Remark}
\newtheorem{ex}[thm]{Example}

\newtheorem{ques}{Question}
\newtheorem*{conj}{Conjecture}

\newtheorem*{ack}{Acknowledgments}

\begin{document}
\title{Relatively spectral morphisms and applications to K-theory}
\author{Bogdan Nica}
\address{\newline Department of Mathematics\newline Vanderbilt University\newline Nashville, TN 37240, USA}
\date{May 30, 2008}

\begin{abstract}
\noin Spectral morphisms between Banach algebras are useful for comparing their K-theory and their ``noncommutative dimensions'' as expressed by various notions of stable ranks. In practice, one often encounters situations where the spectral information is only known over a dense subalgebra. We investigate such relatively spectral morphisms. We prove a relative version of the Density Theorem regarding isomorphism in K-theory. We also solve Swan's problem for the connected stable rank, in fact for an entire hierarchy of higher connected stable ranks that we introduce.
\end{abstract}
\thanks{Supported by FQRNT (Le Fonds qu\'eb\'ecois de la recherche sur la nature et les technologies).}

\maketitle

\section{Introduction}
The following useful criterion for K-theoretic isomorphism is known as the Density Theorem. Initial versions are due to Karoubi \cite[p.109]{Kar78} and Swan \cite[Sec.2.2 \& 3.1]{Swa77}; see also Connes \cite[Appendix 3]{Con81}. The Density Theorem as stated below is taken from Bost \cite[Thm.1.3.1]{Bos90}.

\begin{thm}
Let $\phi:A\To B$ be a dense and spectral morphism between Banach algebras. Then $\phi$ induces an isomorphism $K_*(A)\simeq K_*(B)$.
\end{thm}

\noin A morphism $\phi:A\To B$ is \emph{dense} if $\phi$ has dense image, and \emph{spectral} if $\spec_B(\phi(a))=\spec_A(a)$ for all $a\in A$. Throughout this paper, algebras and their (continuous) morphisms are assumed to be unital.

While proving his version of the density theorem, Swan remarked \cite[p.206]{Swa77} on the possibility that, under the same hypotheses as in the Density Theorem, one has not only isomorphism in K-theory but also equality of stable ranks. That is, \emph{Swan's problem} asks the following:

\begin{ques}  If $\phi:A\To B$ is a dense and spectral morphism between Banach algebras, are the stable ranks of $A$ and $B$ equal?
\end{ques}

There are many notions of stable rank in the literature, and the stable rank in the above problem should be interpreted in a generic sense. Broadly speaking, stable ranks are noncommutative notions of dimension and they are related to stabilization phenomena in K-theory. The stable ranks that are readily interpreted as noncommutative dimensions are the Bass stable rank ($\bsr$) and the topological stable rank ($\tsr$). As for stabilization in K-theory, the most natural rank to consider is the connected stable rank ($\csr$). We thus view Swan's problem for the connected stable rank as the suitable companion to the K-theoretic isomorphism described by the Density Theorem.

Partial results on Swan's problem were obtained by Badea \cite{Bad98}. Most significantly, $\bsr(A)=\bsr(B)$ whenever $A$ is a ``smooth'' subalgebra of a $\Cstar$-algebra $B$ \cite[Thm.1.1, Cor.4.10]{Bad98}. Also, $\csr(A)=\csr(B)$ for $A$ a dense and spectral subalgebra of a commutative Banach algebra $B$ \cite[Thm.4.15]{Bad98}. As for the results concerning the topological stable rank, \cite[Thm.4.13, Cor.4.14]{Bad98},  the hypotheses are unnatural.

In this paper we investigate a weaker notion of spectral morphism.  A morphism $\phi:A\to B$ is \emph{relatively spectral} if $\spec_B(\phi(x))=\spec_A(x)$ for all $x$ in some dense subalgebra $X$ of $A$. We do not know examples of relatively spectral morphisms that are not spectral. Most likely, examples exist and they are not obvious. If $A$ enjoys some form of spectral continuity then a relatively spectral morphism $\phi:A\To B$ is in fact spectral (Section~\ref{speccont}).

The point of considering relatively spectral morphisms is that one can get by with less spectral information. For instance, when one compares two completions of a group algebra $\C\G$, it suffices to consider the spectral behavior of finitely-supported elements. The generalization of the Density Theorem to the relatively spectral context reads as follows:

\begin{thm} Let $\phi:A\To B$ be a dense and completely relatively spectral morphism between Banach algebras. Then $\phi$ induces an isomorphism $K_*(A)\simeq K_*(B)$.
\end{thm}

Following a standard practice, we call a morphism $\phi:A\To B$ \emph{completely} relatively spectral if each amplified morphism $\M_n(\phi): \M_n(A)\To\M_n(B)$ is relatively spectral. It is known (\cite[Lem.2.1]{Swa77},\cite[Prop.A.2.2]{Bos90}, \cite[Thm.2.1]{Sch92}) that spectrality is preserved under amplifications. More precisely, if $\phi:A\To B$ is a dense and spectral morphism then each $\M_n(\phi): \M_n(A)\To\M_n(B)$ is a dense and spectral morphism. We have been unable to prove a similar result for relatively spectral morphisms. Nevertheless, in concrete situations one often encounters a strong form of relative spectrality which propagates to all matrix levels; see Section~\ref{subexp}. 

The surjectivity part in the Relative Density Theorem was known to Lafforgue (\cite[Lem.3.1.1]{Laf07} and comments after \cite[Cor.0.0.3]{Laf02}).

Next we consider Swan's problem for the connected stable rank in the context of relatively spectral morphisms. The answer seems to be the most satisfactory result on Swan's problem so far: 

\begin{thm} Let $\phi:A\To B$ be a dense and relatively spectral morphism between Banach algebras. Then $\csr (A)=\csr(B)$.
\end{thm}

In fact, we prove more. Homotopy stabilization phenomena that are intimately connected to stabilization in K-theory suggest the consideration of higher analogues of connected stable ranks. We show that a relatively spectral, dense morphism preserves these higher connected stable ranks; see Proposition~\ref{csr}.

The Relative Density Theorem can also be considered in a more general context. We introduce certain spectral K-functors $K_\Omega$, indexed by open subsets $\Omega\subseteq\C$ containing the origin. For suitable $\Omega$, one recovers the usual $K_0$ and $K_1$ functors. We prove the Relative Density Theorem for these spectral K-functors; see Proposition~\ref{Kisom}. A quick proof for the usual K-theory is given in Proposition~\ref{quick}.

Handling Swan's problem for higher connected stable ranks and proving the Relative Density Theorem in spectral K-theory is in fact elementary and can be made rather short. The slogan is that dense, relatively spectral morphisms behave well with respect to homotopy, for homotopy of open subsets is equivalent to piecewise-affine homotopy over a dense subalgebra.

Most effort in the sections devoted to the spectral K-functors and the higher connected stable ranks goes towards providing a context for these notions and investing them with meaning. The basic motivation for introducing spectral K-functors is the need for an alternate picture of $K_0$, in which the subset of idempotents is replaced by an open subset. As for the higher connected stable ranks, they are meant to substantiate our claim that, from the K-theoretic perspective, the connected stable rank is the natural stable rank to be considered in Swan's problem. We believe that these two notions, spectral K-functors and higher connected stable ranks, are of independent interest.

Although we are mainly interested in Banach algebras, the results mentioned above are actually true for good Fr\'echet algebras. Section~\ref{goodfrechet} contains the basic facts on good Fr\'echet algebras that are used in this paper. In Section~\ref{RSM} we discuss relatively spectral morphisms. In Section~\ref{finiteness} we show that dense, relatively spectral morphisms transfer the property of being a finite algebra. Section~\ref{homotopy} provides the key homotopy lemma used in proving the Relative Density Theorem and in settling Swan's problem for the higher connected stable ranks. We discuss the higher connected stable ranks in Section~\ref{higherconnected}, and the spectral K-functors in Section~\ref{spectralK}. We close with some applications in Section~\ref{App}. 

\begin{ack} I would like to thank Guoliang Yu for useful discussions, and J\'{a}n \v{S}pakula for a careful reading of the paper. I am grateful to Qayum Khan for explaining me Lemma~\ref{khan}. The fact that Lemma~\ref{khan} holds for compact metric spaces was pointed out to me by Bruce Hughes. Finally, I thank the referee for noticing a faulty argument in the original version.
\end{ack}


\section{Good Fr\'echet algebras}\label{goodfrechet}
Fr\'echet algebras appear naturally in Noncommutative Geometry \cite{Con94}. A motivating example is $C^\infty(M)$, the Fr\'echet algebra of smooth functions on a compact manifold $M$. The context of good Fr\'echet algebras, a context more general than that of Banach algebras, turns out to be the most convenient for our purposes.

\subsection{Good topological algebras} The definition below uses the terminology of \cite[Appendix]{Bos90}:
\begin{defn}\label{good}
A topological algebra $A$ is a \emph{good topological algebra} if the group of invertibles $A^\times$ is open, and the inversion $a\mapsto a^{-1}$ is continuous on $A^\times$. 
\end{defn}

There are two key facts about good topological algebras that we need in what follows:

\begin{prop}\label{goodmatrix}
If $A$ is a good topological algebra then $\M_n(A)$ is a good topological algebra.
\end{prop}

\begin{prop}\label{goodspec} Let $A$ be a good topological algebra. Then:

i) $\spec(a)$ is compact for all $a\in A$;

ii) $A_\Omega=\{a:\spec(a)\subseteq \Omega\}$ is open in $A$ whenever $\Omega\subseteq\C$ is an open set.
\end{prop}

Proposition~\ref{goodmatrix} is due to Swan \cite[Cor.1.2]{Swa77}. The proof of Proposition~\ref{goodspec} is easy and we omit it. 

\subsection{Fr\'echet algebras}
In this paper, we adopt the following:

\begin{defn}\label{frechet} An algebra $A$ is a \emph{Fr\'echet algebra} if $A$ is equipped with a countable family of submultiplicative seminorms $\{\|\cdot\|_k\}_{k\geq 0}$  which make $A$ into a Fr\'echet space.
\end{defn}

We do not require the seminorms to be unital. One can assume, without loss of generality, that the countable family of seminorms in the previous definition is increasing.

If $A$ is a Fr\'echet algebra under the seminorms $\{\|\cdot\|_k\}_{k\geq 0}$ then $\M_n(A)$ is a Fr\'echet algebra under the seminorms given by $\|(a_{ij})\|_k=\sum_{i,j}\|a_{ij}\|_k$; this will be the standard Fr\'echet structure on matrix algebras in what follows.

Every Fr\'echet algebra can be realized as an inverse limit of Banach algebras; this is the Arens - Michael theorem. Indeed, let $A$ be a Fr\'echet algebra under the seminorms $\{\|\cdot\|_k\}_{k\geq 0}$. Let $A_k$ be the Banach algebra obtained by completing $A$ modulo the vanishing ideal of $\|\cdot\|_k$. We obtain an inverse system of Banach algebras with dense connecting morphisms. Then $A$ is isometrically isomorphic to the inverse limit $\varprojlim A_k\subseteq \prod A_k$. Conversely, an inverse limit of Banach algebras is a Fr\'echet algebra: the product $\prod A_k$ of the Banach algebras $A_k$ has a natural Fr\'echet algebra structure given by the coordinate norms, which Fr\'echet structure is inherited by the closed subalgebra $\varprojlim A_k$. 

It is easily checked that $\varprojlim A_k$ is a spectral subalgebra of $\prod A_k$, i.e., $(a_k)\in \varprojlim A_k$ is invertible in  $\varprojlim A_k$ if and only if each $a_k$ is invertible in $A_k$. Thus $\spec\big((a_k)\big)=\cup_k\; \spec_{A_k}(a_k)$ and $r\big((a_k)\big)=\sup_k\; r_{A_k}(a_k)$ for $(a_k)\in \varprojlim A_k$. In particular, if $A$ is a Fr\'echet algebra then $\spec(a)$ is nonempty for all $a\in A$. We also have the following

\begin{prop}
Let $A$ be a Fr\'echet algebra. If $r_A(1-a)<1$ then $a\in A^\times$.
\end{prop}

Viewing a Fr\'echet algebra $A$ once again as an inverse limit of Banach algebras, it is apparent that inversion is continuous on $A^\times$. However, $A^\times$ may not be open. A simple example is $C(\R)$ with the Fr\'echet structure given by the seminorms $\|f\|_{k}=\sup_{x\in [-k,k]} |f(x)|$. The invertible group of $C(\R)$, consisting of the non-vanishing continuous functions, is not open: if $f_k$ is a continuous function such that $f_k=1$ on $[-k,k]$ and $f_k=0$ outside $[-(k+1), k+1]$, then $(f_k)$ is a sequence of non-invertibles converging to $1$.

Consider, on the other hand, the Fr\'echet algebra $C^\infty(M)$ of smooth functions on a compact manifold $M$. The Fr\'echet structure on $C^\infty(M)$ is given by the norms $\|f\|_k=\sum_{|\alpha |\leq k}\|\partial^\alpha f\|_{\infty}$, defined using local charts on $M$. That $C^\infty(M)$ has an open group of invertibles follows by viewing $C^\infty(M)$ as a spectral, continuously-embedded subalgebra of the $\Cstar$-algebra $C(M)$.

If $X$, $Y$ are topological spaces then $X(Y)$ denotes the continuous maps from $Y$ to $X$. Let $\Sigma$ be a compact Hausdorff space and let $A$ be a Fr\'echet algebra under the seminorms $\{\|\cdot\|_k\}_{k\geq 0}$. Then $A(\Sigma)$ is a Fr\'echet algebra under the seminorms $\|f\|_k:=\sup_{p\in\Sigma}\|f(p)\|_k$. As inversion is continuous in $A^\times$, we have $A(\Sigma)^\times=A^\times(\Sigma)$. Since $V(\Sigma)$ is open in $A(\Sigma)$ whenever $V$ is an open subset of $A$, we obtain in particular that $A(\Sigma)$ is good whenever $A$ is good.

\begin{prop}\label{tag}
Let $A$ be a good Fr\'echet algebra, $a\in A$ and $\Omega\subseteq\C$ an open neighborhood of $\spec(a)$. Then there is a unique morphism $\hol(\Omega)\To A$ sending $\id_\Omega$ to $a$, given by
\[\hol_a(f)=f(a):=\frac{1}{2\pi i}\oint f(\lambda)(\lambda-a)^{-1}d\lambda\]
where the integral is taken around a cycle (finite union of closed paths) in $\Omega$ containing $\spec (a)$ in its interior. Furthermore, we have $\spec f(a)=f(\spec (a))$ for each $f\in\hol(\Omega)$.
\end{prop}

\noin The unique morphism indicated by the previous proposition is referred to as the holomorphic calculus for $a$. Here $\hol(\Omega)$, the unital algebra of functions that are holomorphic in $\Omega$, is endowed with the topology of uniform convergence on compacts.


\section{Relatively spectral morphisms}\label{RSM}
In this section, we discuss relatively spectral morphisms. The emphasis is on the comparison between relatively spectral morphisms and spectral morphisms.

Recall that a morphism $\phi: A\To B$ is \emph{spectral} if $\spec_B(\phi(a))=\spec_A(a)$ for all $a\in A$; equivalently, for $a\in A$ we have $a\in A^\times\Leftrightarrow \phi(a)\in B^\times$. We are concerned with the following relative notion:

\begin{defn}
A morphism $\phi: A\To B$ is \emph{spectral relative to a subalgebra} $X\subseteq A$ if $\spec_B(\phi(x))=\spec_A(x)$ for all $x\in X$; equivalently, for $x\in X$ we have $x\in A^\times\Leftrightarrow \phi(x)\in B^\times$. A morphism $\phi: A\To B$ is \emph{relatively spectral} if $\phi$ is spectral relative to some dense  subalgebra of $A$.
\end{defn}

If $\phi: A\To B$ and $\psi: B\To C$ are morphisms, then $\psi\phi$ is spectral relative to $X$ if and only if $\phi$ is spectral relative to $X$ and $\psi$ is spectral relative to $\phi(X)$. This shows, in particular, that the passage from surjective to dense morphisms naturally entails a passage from spectral to relatively spectral morphisms. It also follows that a morphism $\phi: A\To B$ is relatively spectral if and only if both the dense morphism $\phi: A \to \overline{\phi(A)}$ and the inclusion $\overline{\phi(A)}\into B$ are relatively spectral (where $\overline{\phi(A)}$ is the closure of $\phi(A)$ in $B$). In other words, relative spectrality has two aspects: the dense morphism case and the closed-subalgebra case. We are interested in the dense morphism case in this paper.

In practice, the following criterion for relative spectrality is useful:
\begin{prop}\label{radius}
Let $\phi:A \To B$ be a dense morphism between good Fr\'echet algebras. Let $X\subseteq A$ be a dense subalgebra. The following are equivalent:

i) $\spec_B(\phi(x))=\spec_A(x)$ for all $x\in X$

ii) $r_B(\phi(x))=r_A (x)$ for all $x\in X$.
\end{prop}

\begin{proof} i) $\Rightarrow$ ii) is trivial. For ii) $\Rightarrow$ i), let $x\in X$ with $\phi(x)\in B^\times$. Let $(x_n)\subseteq X$ with $\phi(x_n)\To \phi(x)^{-1}$. We have $r_B(1-\phi(x_n)\phi(x))\To 0$, that is, $r_A(1-x_nx)\To 0$. In particular, $x_nx\in A^\times$ for large $n$ so $x$ is left-invertible. A similar argument shows that $x$ is right-invertible. Thus $x\in A^\times$. \end{proof}

Typically, the domain of a dense, relatively spectral morphism cannot be a $\Cstar$-algebra:

\begin{lem}\label{isom} Let $\phi:A\to B$ be a  $*$-morphism, where $A$ is a $\Cstar$-algebra and $B$ is a Banach $*$-algebra. If $\phi$ is dense and spectral relative to a dense $*$-subalgebra, then $\phi$ is an isomorphism.
\end{lem}

\begin{proof} Note first that $\phi$ is onto, as $\phi(A)$ is both dense and closed. Let $X$ be a dense $*$-subalgebra of $A$ relative to which $\phi$ is spectral. For $x\in X$ we have:
\[\|x\|_A^2=\|xx^*\|_A=r_A(xx^*)=r_B(\phi(xx^*))\leq \|\phi(xx^*)\|_B\leq \|\phi(x)\|_B^2\]
It follows that $\|a\|_A\leq \|\phi(a)\|_B\leq C \|a\|_A$ for all $a\in A$, where $C>0$. Thus $\phi$ is an algebraic isomorphism, and $\phi$ can be made into an isometric isomorphism by re-norming $B$.
\end{proof}

Let us point out two disadvantages in working with relatively spectral morphisms.

First, if $\phi:A \To B$ is a spectral morphism between Fr\'echet algebras and $B$ is good, then $A$ is good as well. If $\phi$ is only relatively spectral, then the knowledge that $A$ is good has to come from elsewhere; that is why our applications in Section~\ref{App} involve Banach algebras only. 

Second, spectrality is well-behaved under amplifications: if $\phi:A\To B$ is a dense and spectral morphism then each $\M_n(\phi): \M_n(A)\To\M_n(B)$ is a dense and spectral morphism (\cite[Lem.2.1]{Swa77},\cite[Prop.A.2.2]{Bos90}, \cite[Thm.2.1]{Sch92}). We have been unable to prove a similar result for relatively spectral morphisms.

\begin{ques} 
Let $\phi:A\To B$ be a dense and relatively spectral morphism. Is $\M_n(\phi): \M_n(A)\To\M_n(B)$ a relatively spectral morphism for each $n\geq 1$?
\end{ques}

We thus have to introduce a stronger property that describes relative spectrality at all matrix levels:

\begin{defn}
A morphism $\phi: A\To B$ is \emph{completely spectral relative to a subalgebra} $X\subseteq A$ if each $\M_n(\phi):\M_n(A)\To\M_n(B)$ is spectral relative to $\M_n(X)$. A morphism $\phi: A\To B$ is \emph{completely relatively spectral} if each $\M_n(\phi):\M_n(A)\To\M_n(B)$ is relatively spectral.
\end{defn}

\begin{ex} Consider the dense inclusion $\ell^1\G\into\Cred\G$ for a finitely-generated amenable group $\G$.

If $\G$ has polynomial growth then $\ell^1\G\into\Cred\G$ is spectral (Ludwig \cite{Lud79}). On the other hand, if $\G$ contains a free subsemigroup on two generators then $\ell^1\G\into\Cred\G$ is not spectral (Jenkins \cite{Jen70}). In between these two results, say for groups of intermediate growth, it is unknown whether $\ell^1\G\into\Cred\G$ is spectral or not.

Turning to relative spectrality, it is easy to see that $\ell^1\G\into\Cred \G$ is spectral relative to $\C\G$ if $\G$ has subexponential growth. Indeed, we show that $r_{\ell^1\G}(a)=r_{\Cred\G}(a)$ for $a\in\C\G$. We have $r_{\Cred\G}(a)\leq r_{\ell^1\G}(a)$ from $\|a\|\leq \|a\|_1$, so it suffices to prove that $r_{\ell^1\G}(a)\leq r_{\Cred\G}(a)$. If $a\in\C\G$ is supported in a ball of radius $R$, then $a^n$ is supported in a ball of radius $nR$. We then have 
\[\|a^n\|_1\leq \sqrt{\mathrm{vol}\;B(nR)}\|a^n\|_2\leq  \sqrt{\mathrm{vol}\;B(nR)}\|a^n\|.\] Taking the $n$-th root and letting $n\to\infty$, we obtain $r_{\ell^1\G}(a)\leq r_{\Cred\G}(a)$.

This example serves as a preview for Example~\ref{subgrowth}, where we show that $\ell^1\G\into\Cred \G$ is in fact completely spectral relative to $\C\G$, and for Section~\ref{Sigma1}, where we investigate the groups $\G$ for which $\ell^1\G\into\Cred \G$ is spectral relative to $\C\G$.
\end{ex}

In general, it is hard to decide whether a relatively spectral morphism is spectral or not. For instance, we do not know any examples of relatively spectral morphisms that are not spectral. Under spectral continuity assumptions, however, it is easy to show that relatively spectral morphisms are in fact spectral. We discuss this point in the next subsection.

\subsection{Spectral continuity}\label{speccont} A relatively spectral morphism $\phi: A\To B$ is described by a spectral condition over a dense subalgebra of $A$. In the presence of spectral continuity, this spectral condition can be then extended to the whole of $A$, i.e., $\phi$ is a spectral morphism. It might seem, at first sight, that both $A$ and $B$ need to have spectral continuity for this to work, but in fact spectral continuity for $A$ suffices. 

Spectral continuity can be interpreted in three different ways. The strongest form is to view the spectrum as a map from a good Fr\'echet algebra to the non-empty compact subsets of the complex plane, and to require continuity of the spectrum with respect to the Hausdorff distance. Recall, the Hausdorff distance between two non-empty compact subsets of the complex plane is given by $d_H(C, C')=\inf\{\varepsilon > 0 :\; C\subseteq C'_{\varepsilon}, C'\subseteq C_{\varepsilon}\}$, where $C_{\varepsilon}$ denotes the open $\varepsilon$-neighborhood of $C$. Letting $\mathcal{S}$ denote the class of good Fr\'echet algebras with continuous spectrum in this sense, we have:

\begin{prop}\label{S} Let $\phi:A\To B$ be a relatively spectral morphism between good Fr\'echet algebras. If $A\in\mathcal{S}$ then $\phi$ is spectral.
\end{prop}

\begin{proof} Let $a\in A$; we need to show that $\spec_A(a)\subseteq\spec_B(\phi(a))$. Let $\varepsilon>0$ and pick $x_n\To a$ such that $\spec_A(x_n)=\spec_B(\phi(x_n))$. For $n\gg 1$ we have $\spec_B(\phi(x_n))\subseteq \spec_B(\phi(a))_{\varepsilon}$ by Proposition~\ref{goodspec} ii). On the other hand, the continuity of the spectrum in $A$ gives $\spec_A(a)\subseteq \spec_A(x_n)_{\varepsilon}$ for $n\gg 1$. Combining these two facts, we get $\spec_A(a)\subseteq\spec_A(x_n)_{\varepsilon}=\spec_B(\phi(x_n))_{\varepsilon}\subseteq \spec_B(\phi(a))_{2\varepsilon}$ for $n\gg 1$. As $\varepsilon$ is arbitrary, we are done. \end{proof}

The usefulness of Proposition~\ref{S} is limited by the knowledge about the class $\mathcal{S}$. There are surprisingly few results on $\mathcal{S}$; we did not find in the literature any results complementing the ones contained in Aupetit's survey \cite{Aup79}. The following list summarizes the results mentioned by Aupetit:

a) if $A$ is a commutative Banach algebra then $A$ is in $\mathcal{S}$

b) if $A$ is a commutative Banach algebra and $B$ is a Banach algebra in $\mathcal{S}$, then the projective tensor product $A\otimes B$ is in $\mathcal{S}$ (\cite[Thm.5, p.139]{Aup79})

c) if $A$ is a commutative Banach algebra then $\M_n(A)\in\mathcal{S}$  (\cite[Cor.1, p.139]{Aup79})

d) if $A$ is a Banach algebra in $\mathcal{S}$, then every closed subalgebra of $A$ is in $\mathcal{S}$ (\cite[Thm.3, p.138]{Aup79})

A result of Kakutani says that $B(H)$, the algebra of bounded operators on an infinite-dimensional Hilbert space, is not in $\mathcal{S}$ (\cite[p.34]{Aup79}). 

The second form of spectral continuity, weaker than the one above, is adapted to $*$-algebras. We now require continuity of the spectrum (with respect to the Hausdorff distance) on self-adjoint elements only. Denoting by $\mathcal{S}_\mathrm{sa}$ the class of good Fr\'echet $*$-algebras whose spectrum is continuous on the self-adjoint elements, we have:

\begin{prop} Let $\phi:A\To B$ be a relatively spectral $*$-morphism between good Fr\'echet $*$-algebras. If $A\in\mathcal{S}_\mathrm{sa}$ then $\phi$ is spectral.
\end{prop}

\begin{proof} Arguing as in the proof of Proposition~\ref{S}, we see that $\spec_A(a^*a)=\spec_B(\phi(a^*a))$ for all $a\in A$. Let $a\in A$ with $\phi(a)\in B^\times$. Then $\phi(a^*a), \phi(aa^*)\in B^\times$, and the previous equality of spectra gives $a^*a, aa^*\in A^\times$. Therefore $a\in A^\times$.
 \end{proof}
 
 According to \cite[Cor.4, p.143]{Aup79}, symmetric Banach $*$-algebras are in $\mathcal{S}_\mathrm{sa}$. Recall that a Banach $*$-algebra is \emph{symmetric} if every self-adjoint element has real spectrum. We obtain:
 
\begin{cor} Let $\phi:A\To B$ be a relatively spectral $*$-morphism between Banach $*$-algebras. If $A$ is symmetric then $\phi$ is spectral.
\end{cor}

The third meaning that one can give to spectral continuity is continuity of the spectral radius. Let $\mathcal{R}$ denote the class of good Fr\'echet algebras with continuous spectral radius. Continuity of the spectrum (with respect to the Hausdorff distance) implies continuity of the spectral radius, that is, $\mathcal{S}\subseteq\mathcal{R}$; the inclusion is strict (\cite[p.38]{Aup79}). Kakutani's result, mentioned above, actually says that $B(H)$ is not in $\mathcal{R}$.

\begin{prop}
Let $\phi:A \To B$ be a dense and relatively spectral morphism between good Fr\'echet algebras. If $A\in\mathcal{R}$ then $\phi$ is spectral.
\end{prop}

\begin{proof} The proof is very similar to that of Proposition~\ref{S}. Let $a\in A$; we show that $r_A(a)\leq r_B(\phi(a))$. For then Proposition~\ref{radius} allows us to conclude that $\phi$ is spectral. Let $\varepsilon>0$ and pick $x_n\To a$ such that $r_A(x_n)=r_B(\phi(x_n))$. For $n\gg 1$ we have $r_B(\phi(x_n))\leq r_B(\phi(a))+\varepsilon$, i.e., $r_A(x_n)\leq r_B(\phi(a))+\varepsilon$. Letting $n\To\infty$ and using the continuity of the spectral radius in $A$, we get $r_A(a)\leq r_B(\phi(a))+\varepsilon$. As $\varepsilon$ is arbitrary, we are done.
\end{proof}

Finally, note that if $B$ satisfies one of the above forms of spectral continuity and $\phi: A\To B$ is spectral, then $A$ satisfies that spectral continuity as well.


\section{Finiteness}\label{finiteness}
Recall that an algebra $A$ is \emph{finite} if the left-invertibles of $A$ are actually invertible, equivalently, if the right-invertibles of $A$ are actually invertible. This terminology is standard in the Banach- and $\Cstar$-algebraic setting; in noncommutative ring theory, this property is called Dedekind-finiteness.  An algebra $A$ is \emph{stably finite} if each matrix algebra $\M_n(A)$ is finite. Tracial $\Cstar$-algebras are important examples of stably finite algebras. Stable-finiteness is relevant in K-theory, for it guarantees the non-vanishing of $K_0$.

\begin{prop}\label{finite} Let $\phi: A\To B$ be a dense morphism between good topological algebras. 

a) Assume $\phi$ is relatively spectral. Then $A$ is finite if and only if $B$ is finite.

b) Assume $\phi$ is completely relatively spectral. Then $A$ is stably finite if and only if $B$ is stably finite.
\end{prop}

\begin{proof} We prove a). Part b) is an obvious corollary.

Letting $L(A)$ denote the left-invertibles of $A$, note that the density of $A^\times$ in $L(A)$ suffices for finiteness. Indeed, if $a\in L(A)$ then $u_n\To a$ for some invertibles $u_n$. There is $a'\in A$ with $a'a=1$, so $a'u_n\To 1$. Hence $a'u_n\in A^\times$ for large $n$, thus $a'\in A^\times$ and we conclude $a\in A^\times$.

Let $X$ be a dense subalgebra of $A$ relative to which $\phi$ is spectral.

Assume $B$ is finite. Let $a\in L(A)$ and let $x_n\To a$ where $x_n\in X$. Then $\phi(x_n)\To \phi(a)\in L(B)=B^\times$ so, by relative spectrality, $x_n$ is eventually invertible. Thus $A^\times$ is dense in $L(A)$. 

Assume $A$ is finite. We claim that $\phi(L(A))$ is dense in $L(B)$; as $\phi(L(A))=\phi(A^\times)\subseteq B^\times\subseteq L(B)$, it will follow that $B^\times$ is dense in $L(B)$. Let $b\in L(B)$ with $b'b=1$. Pick sequences $(x_n)$, $(x'_n)$ in $X$ so that $\phi(x_n)\To b$, $\phi(x'_n)\To b'$. Since $\phi(x'_nx_n)\To 1$, relative spectrality gives that $x_n$ is left-invertible for large $n$. This proves that $\phi(L(A))$ is dense in $L(B)$.
\end{proof}


\section{A homotopy lemma}\label{homotopy}
We adapt Proposition A.2.6 of \cite{Bos90} as follows:

\begin{lem}\label{key} Let $A$, $B$ be Fr\'echet spaces and $\phi:A \To B$ be a continuous linear map with dense image. Let $U\subseteq A$, $V\subseteq B$ be open with $\phi(U)\subseteq V$. If $X$ is a dense subspace of $A$ such that $U\cap X=\phi^{-1}(V)\cap X$ then $\phi$ induces a weak homotopy equivalence between $U$ and $V$.
\end{lem}

\begin{proof} First, we show that $\phi$ induces a bijection $\pi_0(U)\To\pi_0(V)$. We make repeated use of the local convexity. For surjectivity, pick $v\in V$. As $V$ is open and $\phi(X)$ is dense, there is $v_X\in V\cap\phi(X)$ such that $v$ and $v_X$ can be connected by a segment in $V$. If $u_X\in X$ is a pre-image of $v_X$, then $u_X\in U$. For injectivity, let $u, u'\in U$ such that $\phi(u),\phi(u')$ are connected in $V$. As $U$ is open and $X$ is dense, there are $u_X, u'_X\in U\cap X$ such that $u$ and $u_X$, respectively $u'$ and $u'_X$, can be connected by a segment in $U$. That is, we may assume that we start with $u, u'\in U\cap X$. As $V$ is open and $\phi(X)$ is dense, if $\phi(u)$ and $\phi(u')$ can be connected by a path in $V$ then they can be connected by a piecewise-linear path $p_B$ lying entirely in $V\cap\phi(X)$. Take pre-images in $X$ for the vertices of $p_B$ and extend to a piecewise-linear path $p_A$ connecting $u$ to $u'$. Since the path $p_A$ lies in $X$ and is mapped inside $V$, it follows that $p_A$ lies in $U$. We conclude that $u,u'$ are connected in $U$.

Next, let $k\geq 1$. We show that $\phi$ induces a bijection $\pi_k(U,a)\To\pi_k(V,\phi(a))$ for each $a\in U$. Up to translating $U$ by $-a$ and $V$ by $-\phi(a)$, we may assume that $a=0$. Fix a basepoint $\bullet$ on $S^k$. The dense linear map $\phi:A\To B$ induces a dense linear map $\phi_k: A(S^k)\To B(S^k)$, which restricts to a dense linear map $\phi^\bullet_k: A(S^k)^\bullet\To B(S^k)^\bullet$. The $\bullet$-decoration stands for restricting to the based maps sending $\bullet$ to $0$. With $U(S^k)^\bullet$, $V(S^k)^\bullet$, $X(S^k)^\bullet$ playing the roles of $U$, $V$, $X$, we get by the first part that $\phi_k^\bullet$ induces a bijection $\pi_0(U(S^k)^\bullet)\To \pi_0(V(S^k)^\bullet)$. That means precisely that $\phi$ induces a bijection $\pi_k(U,0)\To\pi_k(V,0)$.

Let us give some details on the previous paragraph. That $A(S^k)^\bullet$, respectively $B(S^k)^\bullet$, is a Fr\'echet space follows by viewing it as a closed subspace of the Fr\'echet space $A(S^k)$, respectively $B(S^k)$. The map $\phi_k: A(S^k)\To B(S^k)$ is given by $\phi_k(f)=\phi\circ f$, and it is obviously linear and continuous. Consequently, the restriction $\phi^\bullet_k: A(S^k)^\bullet\To B(S^k)^\bullet$ is linear and continuous. The fact that $X(S^k)^\bullet$ is dense in $A(S^k)^\bullet$ is showed using a partition of unity argument. This also proves the density of $\phi_k$ and $\phi_k^\bullet$. Next, $U(S^k)$ is open in $A(S^k)$ (see comments before Proposition~\ref{tag}) so $U(S^k)^\bullet$ is open in $A(S^k)^\bullet$; similarly, $V(S^k)^\bullet$ is open in $B(S^k)^\bullet$. Finally, $\phi_k^\bullet \big(U(S^k)^\bullet\big)\subseteq V(S^k)^\bullet$ and $U(S^k)^\bullet\cap X(S^k)^\bullet=(\phi_k^\bullet)^{-1}\big(V(S^k)^\bullet\big)\cap X(S^k)^\bullet$ are immediate to check, they boil down to $\phi(U)\subseteq V$ and $U\cap X=\phi^{-1}(V)\cap X$ respectively.
\end{proof}

It is not hard to imagine that, modulo notational complications, the idea used to treat $\pi_0$ has a higher-dimensional analogue. The underlying phenomenon is that homotopical considerations for open subsets of Fr\'echet spaces can be carried out in a piecewise-affine fashion over a dense subspace. Bost's elegant approach of upgrading $\pi_0$ knowledge to higher homotopy groups is, however, more economical.

As we shall see, the above Lemma immediately yields the Relative Density Theorem in spectral K-theory (Proposition~\ref{Kisom}) and a positive answer to Swan's problem for the higher connected stable ranks (Proposition~\ref{csr}). We also get a quick proof for the Relative Density Theorem in the Banach algebra case:

\begin{cor}\label{quick}
Let $\phi:A\To B$ be a dense and completely relatively spectral morphism between Banach algebras. Then $\phi$ induces an isomorphism $K_*(A)\simeq K_*(B)$.
\end{cor}

\begin{proof}
Let $\phi:A\To B$ be a dense morphism that is spectral relative to a dense subalgebra $X\subseteq A$, that is, $A^\times\cap X=\phi^{-1}(B^\times)\cap X$. Lemma~\ref{key} says that $\phi$ induces a bijection $\pi_*(A^\times)\To\pi_*(B^\times)$ for $*=0,1$. Thus, if $\phi:A\To B$ is a dense and completely relatively spectral morphism then $\phi$ induces a bijection $\pi_*(\GL_n(A))\To\pi_*(\GL_n(B))$. It follows that $\phi$ induces a bijection $\varinjlim\;\pi_*(\GL_n(A))\To\varinjlim\;\pi_*(\GL_n(B))$. That is, $\phi$ induces an isomorphism $K_{*+1}(A)\simeq K_{*+1}(B)$. As the Bott isomorphism $K_0(A)\simeq K_2(A)$ is natural, we conclude that $\phi$ induces an isomorphism $K_*(A)\simeq K_*(B)$.\end{proof}


\section{Higher connected stable ranks}\label{higherconnected}
This section is devoted to higher analogues of the notion of connected stable rank. We start by defining the higher connected stable ranks, and by estimating them in terms of the topological stable rank. We compute the higher connected stable ranks of $C(\Sigma)$ for certain finite CW-complexes $\Sigma$. Then we explain the relation between the higher connected stable ranks and certain homotopy stabilization ranks. In particular, we obtain information about stabilization in K-theory. Finally, we give a positive answer to Swan's problem for the higher connected stable ranks.

\subsection{Defining and estimating higher connected stable ranks}
Stable ranks for a topological algebra $A$ are usually defined in terms of the sequence of left-generating sets
\[\mathrm{Lg}_n(A)=\{(a_1,\dots,a_n):\; Aa_1+\dots + Aa_n=A\}\subseteq A^n\]
where $n\geq 1$. It is easy to see that, for a good topological algebra $A$, $\Lg_n(A)$ is open in $A^n$.

Recall:

\begin{defn} Let $A$ be a good topological algebra.

The \emph{topological stable rank} $\tsr(A)$ is the least $n$ such that $\Lg_n(A)$ is dense in $A^n$.

The \emph{connected stable rank} $\csr(A)$ is the least $n$ such that $\Lg_m(A)$ is connected for all $m\geq n$.

The \emph{Bass stable rank} $\bsr(A)$ is the least $n$ such that $\Lg_{n+1}(A)$ is reducible to $\Lg_n(A)$ in the following sense: if $(a_1, \dots, a_{n+1})\in\Lg_{n+1}(A)$ then $(a_1+x_1a_{n+1},\dots,a_n+x_na_{n+1})\in\Lg_n(A)$ for some $x_1,\dots, x_n$ in $A$.
\end{defn}

\begin{rem} The Bass stable rank, a purely algebraic notion, is due to Bass \cite{Bas64}. The topological stable rank and the connected stable rank were introduced by Rieffel \cite{Rie83}; implicitly, these two stable ranks also appear in \cite{CL84}. Both \cite{Rie83} and \cite{CL84} address the basic inequality  $\csr-1\leq \bsr\leq\tsr$. Note that the definition of the connected stable rank that we adopt in this paper is not the original definition \cite[Def.4.7]{Rie83} but rather an equivalent one \cite[Cor.8.5]{Rie83}. Note also that the stable ranks defined above are, strictly speaking, left stable ranks. As usual, a stable rank is declared to be infinite if there is no integer fulfilling the requirement.
\end{rem}

We define higher connected stable ranks, that is, stable ranks encoding the higher connectivity of the left-generating sets. As we shall see in the next subsection, such higher connectivity properties are relevant for stabilization phenomena in K-theory.

\begin{defn} Let $A$ be a good topological algebra and $k\geq 0$. The \emph{$k$-connected stable rank} $\csr_k(A)$ is the least $n$ such that $\Lg_m(A)$ is $k$-connected for all $m\geq n$.
\end{defn}

Recall, a space $\Sigma$ is said to be $k$-connected if $\pi_i(\Sigma,\bullet)=0$ for all $0\leq i\leq k$ (this is independent of the choice of basepoint $\bullet$ in $\Sigma$). Equivalently, $\Sigma$ is $k$-connected if each map $S^i\To\Sigma$ can be extended to a map $I^{i+1}\To\Sigma$, for all $0\leq i\leq k$. Henceforth, $I$ denotes the unit interval.

This definition yields a hierarchy of connected stable ranks 
\[\csr_0\leq\csr_1\leq \csr_2\leq \dots\]
in which the very first one is the usual connected stable rank; indeed, topological algebras being locally path-connected, connectivity and path-connectivity are equivalent for open subsets. We now seek to generalize the following two important facts concerning the connected stable rank:

i) $\csr$ is homotopy invariant, i.e., $\csr(A)=\csr(B)$ whenever $A$ and $B$ are homotopy equivalent;

ii) $\csr\leq \tsr +1$

Fact i) this is due to Nistor \cite{Nis86}. Recall that two morphisms $\phi_0,\phi_1:A\To B$ are said to be  homotopic if they are the endpoints of a path of morphisms $\{\phi_t\}_{0\leq t\leq 1}:A\To B$; here $t\mapsto \phi_t$ is continuous in the sense that $t\mapsto\phi_t(a)$ is continuous for each $a\in A$. If there are morphisms $\alpha:A\To B$ and $\beta:B\To A$ with $\beta\alpha$ homotopic to $\id_A$ and $\alpha\beta$ homotopic to $\id_B$, then $A$ and $B$ are said to be homotopy equivalent.

Fact ii) can be obtained, for instance, by combining $\csr(A)\leq \tsr (A(I))$ (a fact implicit in \cite{Nis86}) with $\tsr(A(I))\leq \tsr(A)+1$ (\cite[Cor.7.2]{Rie83}). We follow this strategy in the next proposition. Recall that we use $X(Y)$ to denote the continuous maps from $Y$ to $X$.

\begin{prop}\label{csrk}  For each $k\geq 0$ we have:

i) $\csr_k$ is homotopy invariant.

ii) $\csr_k(A)\leq \tsr(A(I^{k+1}))$. In particular, $\csr_k(A)\leq \tsr(A)+k+1$. If Rieffel's estimate
\[(R)\quad\quad\tsr(A(I^2))\leq \tsr(A)+1\]
holds, then in fact $\csr_k(A)\leq \tsr(A)+\big\lfloor\frac{1}{2}k\big\rfloor+1$.
\end{prop}

\begin{proof} i) Homotopic morphisms $\phi_0,\phi_1 :A\To B$ induce homotopic maps $\phi_0,\phi_1:\Lg_m(A)\To \Lg_m(B)$. Thus, if $A$ and $B$ are homotopy equivalent (as algebras), then $\Lg_m(A)$ and $\Lg_m(B)$ are homotopy equivalent (as topological spaces).

ii) We show the slightly better estimate $\csr_k(A)\leq \bsr(A(I^{k+1}))$. To that end we use the following fact \cite[Lem.4.1]{Bas64}: if $A\To B$ is onto and $m\geq\bsr(A)$ then $\Lg_m(A)\To\Lg_m(B)$ is onto. It follows, in particular, that $\bsr(B)\leq \bsr(A)$ whenever $B$ is a quotient of $A$.

Let $m\geq \bsr(A(I^{k+1}))$. We need to show that $\Lg_m(A)$ is $k$-connected. We fix $0\leq i\leq k$ and we show that every map $S^i\To\Lg_m(A)$ extends to a map $I^{i+1}\To\Lg_m(A)$. 

Consider the restriction morphism $A(I^{i+1})\onto A(S^i)$. We have $m\geq \bsr(A(I^{i+1}))$, as $A(I^{i+1})$ is a quotient of $A(I^{k+1})$, so the restriction map $\Lg_m(A(I^{i+1}))\To \Lg_m(A(S^i))$ is onto. For a compact space $\Sigma$ - more generally, for paracompact $\Sigma$ \cite[Lem.3.2]{Bad98} - one has a natural identification $\Lg_m(A(\Sigma))\simeq(\Lg_m(A))(\Sigma)$. We thus have an onto map $(\Lg_m(A))(I^{i+1})\To (\Lg_m(A))(S^i)$, proving the desired extension property.

From $\tsr(A(I))\leq \tsr(A)+1$ one obtains $\tsr(A(I^{k+1}))\leq \tsr(A)+k+1$, hence the first estimate on $\csr_k$ in terms of $\tsr$. If $(R)$ holds, then one obtains $\tsr(A(I^{k+1}))\leq \tsr(A)+\big\lfloor\frac{1}{2}k\big\rfloor+1$, hence the second estimate on $\csr_k$.
\end{proof}

\begin{rem} The estimate $(R)$ appears as Question 1.8 in \cite{Rie83}. The motivation comes from the commutative case, as one has $\tsr(C(\Sigma))=\big\lfloor\frac{1}{2}\dim(\Sigma)\big\rfloor+1$ whenever $\Sigma$ is a compact space. Sudo claims in \cite{Sud03} a proof for $(R)$; however, Proposition 1 of \cite{Sud03} cannot hold. For if it is true, then it implies $(\dagger)\;\; \tsr(A(I))=\max\{\csr(A), \tsr(A)\}$ for all $\Cstar$-algebras $A$, by using Nistor's description for $\tsr(A(I))$ as the absolute connected stable rank of $A$. Putting $A=C(I^k)$ in $(\dagger)$ leads to $\tsr(C(I^{k+1}))=\tsr(C(I^k))$ for all $k\geq 1$, since $\csr (C(I^k))=1$ by the homotopy invariance of the connected stable rank. This is a contradiction. The reader might enjoy finding elementary counterexamples to \cite[Prop.1]{Sud03}.
\end{rem}

\subsection{Higher connected stable ranks for $C(\Sigma)$}\label{compute}  From Proposition~\ref{csrk} we obtain a dimensional upper bound for the higher connected stable ranks: 
\[\csr_k(C(\Sigma))\leq\tsr(C(\Sigma\times I^{k+1}))=\bigg\lfloor\frac{1}{2}(\dim(\Sigma)+k+1)\bigg\rfloor+1\] That is:
\begin{equation}\label{+}
\csr_k(C(\Sigma))\leq\bigg\lceil\frac{\dim(\Sigma)+k}{2}\bigg\rceil+1
\end{equation}
In general, we cannot have equality in \eqref{+}, for the left-hand side is homotopy invariant whereas the right-hand side is not. We will see, however, that in many natural cases we do get equality in \eqref{+}. 

Let us show, for a start, that $\csr_k(\C)=\big\lceil\frac{1}{2}k\big\rceil+1$. As $\Lg_m(\C)=\C^m\setminus\{0\}$ has the homotopy type of $S^{2m-1}$, which is $(2m-2)$-connected but not $(2m-1)$-connected, we obtain that $\csr_k(\C)$ is the least $n$ such that  $2n-2\geq k$, i.e., $\csr_k(\C)=\big\lceil\frac{1}{2}k\big\rceil+1$. It follows that $\csr_k(C(\Sigma))=\big\lceil\frac{1}{2}k\big\rceil+1$ for $\Sigma$ contractible.

The computation of $\csr_k(C(\Sigma))$ is a question about the vanishing of certain cohomotopy groups associated to $\Sigma$. Indeed, if we identify $\Lg_m(C(\Sigma))$ with $(\C^m\setminus\{0\})(\Sigma)$, then for $i\geq 0$ we have:
\[
[S^i, \Lg_m(C(\Sigma))]=[S^i, (\C^m\setminus\{0\})(\Sigma)] =[\Sigma\times S^i, \C^m\setminus\{0\}]=[\Sigma\times S^i, S^{2m-1}]
\]
Thus $\csr_k(C(\Sigma))$ is the least $n$ such that, for all $m\geq n$, we have $[\Sigma\times S^i, S^{2m-1}]=0$ for $0\leq i\leq k$. This explicit formulation shows, once again, that the higher connected stable ranks of $C(\Sigma)$ only depend on the homotopy type of $\Sigma$. We can also recover \eqref{+}. Indeed, if $\Sigma'$ is a finite CW-complex then $[\Sigma', S^N]=0$ as soon as $N>\dim(\Sigma')$. Thus, if $\dim(\Sigma)+k< 2m-1$ then $[\Sigma\times S^i, S^{2m-1}]=0$ for all $0\leq i\leq k$. We obtain $\csr_k(C(\Sigma))\leq\big\lceil\frac{1}{2}(\dim(\Sigma)+k)\big\rceil+1$, i.e., \eqref{+}.

Now, for suitable dimensions, non-trivial cohomology implies non-trivial cohomotopy:

\begin{lem}\label{khan} Let $\Sigma'$ be an $N$-dimensional finite CW-complex. Then:

a) $H^N(\Sigma';\Z)\neq 0$ if and only if $[\Sigma', S^N]\neq 0$;

b) $H^{N-1}(\Sigma';\Z)\neq 0$ implies $[\Sigma', S^{N-1}]\neq 0$.
\end{lem}

\noin Note that part b) is not an equivalence, as witnessed by $\Sigma'=S^N$ for $N\geq 3$. The previous lemma also holds for $\Sigma'$ an $N$-dimensional compact metric space, by results from Hurewicz \& Wallman: \emph{Dimension Theory}, pp. 149-150. The cohomology is then understood in the \v{C}ech sense.

\begin{proof} Recall that $H^N(\Sigma';\Z)$ can be identified with $[\Sigma', K(\Z,N)]$. Realize the $K(\Z,N)$ CW-complex by starting with $S^N$ as its $N$-skeleton and then adding cells of dimension $\geq N+2$ that kill the higher homotopy. 

For a) we claim that the natural map $[\Sigma', S^N]\To [\Sigma', K(\Z,N)]$ is bijective. Surjectivity is clearly a consequence of the Cellular Approximation Theorem. But so is injectivity: a homotopy in $K(\Z,N)$ between $f_0,f_1:\Sigma'\To S^N$ can be assumed cellular, so it maps to the $(N+1)$-skeleton of $K(\Z,N)$ which is $S^N$. That is, the maps $f_0,f_1:\Sigma'\To S^N$ are actually homotopic over $S^N$. 

For b), note that $[\Sigma', S^N]\To [\Sigma', K(\Z,N)]$ remains surjective if $\Sigma'$ is $(N+1)$-dimensional. Again, this follows from the Cellular Approximation Theorem and the fact that the $(N+1)$-skeleton of $K(\Z,N)$ is $S^N$. Replacing $N$ by $N-1$ gives b) as stated. \end{proof}

We obtain the following sufficient condition for having equality in \eqref{+}:

\begin{prop}\label{d} Let $\Sigma$ be a $d$-dimensional finite CW-complex.

a) Let $k\geq 1$. If $H^d(\Sigma;\Z)\neq 0$ then $\csr_k(C(\Sigma))=\big\lceil\frac{1}{2}(d+k)\big\rceil+1$.

b) Let $k=0$. If $H^d(\Sigma;\Z)\neq 0$ for odd $d$, or $H^{d-1}(\Sigma;\Z)\neq 0$ for even $d$, then $\csr(C(\Sigma))=\big\lceil\frac{1}{2}d\big\rceil+1$.
\end{prop}

\begin{proof}
a) Let $m=\big\lceil\frac{1}{2}(d+k)\big\rceil$. Then $\dim(\Sigma\times S^i)=2m-1$ for some $0\leq i\leq k$, namely $i=k-1$ or $i=k$. Now $H^{2m-1}(\Sigma\times S^i;\Z)\simeq H^d(\Sigma;\Z)\otimes_\Z H^i(S^i;\Z)$ is non-vanishing, so $[\Sigma\times S^i, S^{2m-1}]\neq 0$. Thus $\csr_k(C(\Sigma)) > \big\lceil\frac{1}{2}(d+k)\big\rceil$.

b) Let $m=\big\lceil\frac{1}{2}d\big\rceil$. If $d$ is odd, then $2m-1=d$ so $H^d(\Sigma;\Z)\neq 0$ implies $[\Sigma, S^{2m-1}]\neq 0$. Therefore $\csr (C(\Sigma))> \big\lceil\frac{1}{2}d\big\rceil$. If $d$ is even, then $2m-1=d-1$. Here we use the fact that $H^{d-1}(\Sigma;\Z)\neq 0$ implies $[\Sigma, S^{2m-1}]\neq 0$. We get $\csr (C(\Sigma))> \big\lceil\frac{1}{2}d\big\rceil$ in this case as well. 
\end{proof}

\begin{ex}\label{torus} We have $\csr_k(C(T^d))=\big\lceil\frac{1}{2}(d+k)\big\rceil+1$ for the $d$-torus.
\end{ex}

\begin{ex} For the $d$-sphere, Proposition~\ref{d} gives $\csr_k(C(S^d))=\big\lceil\frac{1}{2}(d+k)\big\rceil+1$ for $k\geq 1$, and for $k=0$ and odd $d$. If $k=0$ and $d$ is even, the cohomological criterion fails us. Nevertheless, we have $[S^d,S^{d-1}]\neq 0$ for $d\neq 2$, hence $\csr_k(C(S^d))=\big\lceil\frac{1}{2}(d+k)\big\rceil+1$ in this case, too. Summarizing, we have $\csr_k(C(S^d))=\big\lceil\frac{1}{2}(d+k)\big\rceil+1$ except for $k=0$, $d=2$. We leave it to the reader to verify that $\csr(C(S^2))=1$.
\end{ex}

\subsection{Stabilization in K-theory and higher connected stable ranks}
As before, let $A$ be a good topological algebra. The sequence 
\[\{1\}=\GL_0(A)\into A^\times=\GL_1(A)\into\GL_2(A)\into\dots\] 
induces, for each $i\geq 0$, a sequence of (identity-based) homotopy groups:
\[(\pi_i)\qquad\pi_i(\GL_0(A))\to \pi_i(\GL_1(A))\to\pi_i(\GL_2(A))\to\dots\]
Say that $m$ is a stable level for $(\pi_i)$ if $\pi_i(\GL_{m'}(A))\To\pi_i(\GL_{m'+1}(A))$ is an isomorphism for all $m'\geq m$. In analogy with $k$-connectivity, which requires vanishing of all homotopy groups up to the $k$-th one, the following notion encodes the stabilization of $(\pi_i)$ for all $0\leq i\leq k$:

\begin{defn} Let $k\geq 0$. The \emph{$k$-homotopy stabilization rank} $\hsr_k(A)$ is the least nonnegative integer which is a stable level for $(\pi_i)$, for all $0\leq i\leq k$.
\end{defn}

Homotopy stabilization ranks are closely related to higher connected stable ranks:

\begin{prop} $\hsr_k(A)\leq\csr_{k+1}(A)-1$ and $\csr_k(A)-1\leq\max\{\hsr_k(A), \csr(A)-1\}$.
\end{prop}

Roughly speaking, we have $\csr_k(A)-1\leq \hsr_k(A)\leq\csr_{k+1}(A)-1$. This actually holds if, say, $\tsr(A)=1$ and $K_1(A)\neq 0$, for then $\hsr_0(A)\geq 1\geq \csr(A)-1$.

\begin{proof}
Let $m\geq \csr(A)-1$. According to \cite[Cor.1.6]{CL86} one has a long exact homotopy sequence:
\begin{eqnarray*} 
\dots \to\pi_{i+1}(\Lg_{m+1}(A))\to\pi_i(\GL_m(A))\to\pi_i(\GL_{m+1}(A))\to\pi_i(\Lg_{m+1}(A))\to\cdots\\ \dots\to\pi_1(\Lg_{m+1}(A))\to\pi_0(\GL_m(A))\to\pi_0(\GL_{m+1}(A))\to 0\phantom{xxxxxxxxxxxx}
\end{eqnarray*}
Let $m\geq \csr_{k+1}(A)-1$. Then $\pi_i(\GL_m(A))\to\pi_i(\GL_{m+1}(A))$ is an isomorphism for $0\leq i\leq k$. Hence $\csr_{k+1}(A)-1$ is a stable level for $(\pi_i)$, for all $0\leq i\leq k$, so $\hsr_k(A)\leq \csr_{k+1}(A)-1$. 

For the second inequality, let $m\geq \max\{\hsr_k(A), \csr(A)-1\}$. Then $\Lg_{m+1}(A)$ is connected and, for all $1\leq i\leq k$, $\pi_i(\Lg_{m+1}(A))$ is trivial since it is squeezed between two isomorphisms in the above exact sequence. Thus $\Lg_{m+1}(A)$ is $k$-connected. \end{proof}

\begin{cor} Let $\Sigma$ be a finite CW-complex of dimension $2d$ whose top cohomology group is non-zero. Then $\hsr_1(C(\Sigma))=d+1$, i.e., the level at which both $(\pi_0)$ and $(\pi_1)$ begin to stabilize is $d+1$. 
\end{cor}

\begin{cor} \label{hest} $\hsr_k(A)\leq \tsr(A)+k+1$; if $(R)$ holds, then $\hsr_k(A)\leq \tsr(A)+\big\lceil\frac{1}{2}k\big\rceil$.
\end{cor}

\begin{rem} The estimate $\hsr_k(A)\leq \tsr(A)+k+1$ is essentially Theorem 6.3 of \cite{CL86}. However,  such an estimate is not our main goal. The purpose of the above discussion is rather to emphasize the following ideas:

i) the stable ranks that are best suited for controlling stabilization in the homotopy sequences $(\pi_*)$ are the connected stable ranks;

ii) the connected stable ranks can be estimated by the ``dimensional'' stable ranks, namely the topological stable rank and the Bass stable rank, which are typically easier to compute;

iii) Rieffel's conjectured inequality $(R)$ plays a crucial role in providing good estimates for the connected stable ranks, hence in estimating homotopy stabilization.
\end{rem}

For Banach algebras, the K-theoretic interpretation is the following. The $K_1$ group and the $K_0$ group are the limit groups of the direct sequences $(\pi_0)$ and $(\pi_1)$:
\[K_1(A)=\varinjlim\; \pi_0(\GL_n(A)), \qquad K_0(A)\simeq \varinjlim\; \pi_1(\GL_n(A))\]
Let us stress the following aspect: stabilization for $(\pi_0)$ is indeed stabilization for $K_1$, whereas stabilization for $(\pi_1)$ is construed as stabilization for $K_0$.

Therefore:

$\cdot$ $K_1(A)\simeq \pi_0(\GL_n(A))$ for $n\geq\csr_1(A)-1$;

$\cdot$ $K_1(A)\simeq \pi_0(\GL_n(A))$ and $K_0(A)\simeq \pi_1(\GL_n(A))$ for $n\geq\csr_2(A)-1$.

\noin In particular, relative to the topological stable rank we have

$\cdot$ $K_1(A)\simeq \pi_0(\GL_n(A))$ for $n\geq \tsr(A)+1$ (compare \cite[Thm.10.12]{Rie83});

$\cdot$ $K_1(A)\simeq \pi_0(\GL_n(A))$ and $K_0(A)\simeq \pi_1(\GL_n(A))$ for $n\geq \tsr(A)+2$.

\noin If, furthermore, $(R)$ holds, then

$\cdot$ $K_1(A)\simeq \pi_0(\GL_n(A))$ for $n\geq \tsr(A)$;

$\cdot$ $K_1(A)\simeq \pi_0(\GL_n(A))$ and $K_0(A)\simeq \pi_1(\GL_n(A))$ for $n\geq\tsr(A)+1$.

\subsection{Swan's problem for the higher connected stable ranks}\label{Swan} Recall that Swan's problem asks the following: if $\phi:A \to B$ is a dense and spectral morphism, are the stable ranks of $A$ and $B$ equal? We show that this is indeed the case for the higher connected stable ranks. The next result improves a result of Badea \cite[Thm.4.15]{Bad98} in several ways. First and foremost, it removes the commutativity assumption. Second, it suffices to know that the morphism is relatively spectral. Third, it handles the whole hierarchy of connected stable ranks.

\begin{prop}\label{csr} Let $\phi:A\To B$ be a dense and relatively spectral morphism, where $A$ and $B$ are good Fr\'echet algebras. Then $\csr_k (A)=\csr_k(B)$ for all $k\geq 0$.
\end{prop}

\begin{proof} 
We have $\phi(\Lg_m(A))\subseteq \Lg_m(B)$. We claim that $\Lg_m(A)\cap X^m=\phi^{-1} (\Lg_m(B))\cap X^m$, where $X$ is a dense subalgebra of $A$ relative to which $\phi$ is spectral. Let $(x_i)\in X^m$ with $(\phi(x_i))\in \Lg_m(B)$, so $\sum c_i\phi(x_i)=1$ for some $(c_i)\in B^m$. Density of $\phi$ allows to approximate each $c_i$ by some $\phi(x_i')$, where $x_i'\in X$, so as to get $\sum\phi(x'_i)\phi(x_i)=\phi\big(\sum x'_ix_i\big)\in B^\times$. As $\sum x'_ix_i\in X$, we infer that $\sum x'_ix_i\in A^\times$ by relative spectrality. Thus $(x_i)\in\Lg_m(A)$ as desired.

Lemma~\ref{key} applied to the dense morphism $\phi: A^m\To B^m$ gives that $\phi$ induces a weak homotopy equivalence between $\Lg_m(A)$ and $\Lg_m(B)$. In particular, $\Lg_m(A)$ is $k$-connected if and only if $\Lg_m(B)$ is $k$-connected. \end{proof}

\begin{ex} We have $\csr_k(\Cred(\Z^d))=\big\lceil\frac{1}{2}(d+k)\big\rceil+1$ by Example~\ref{torus}. Since the dense inclusion $\ell^1(\Z^d)\into\Cred(\Z^d)$ is spectral, it follows that $\csr_k(\ell^1(\Z^d))=\big\lceil\frac{1}{2}(d+k)\big\rceil+1$ as well.
\end{ex}


\section{Spectral K-functors}\label{spectralK}
By a \emph{K-functor} we simply mean a functor from good Fr\'echet algebras to abelian groups. The notion of K-scheme we introduce below provides a general framework for constructing K-functors. A K-functor is usually required to be stable, homotopy-invariant, half-exact and continuous. These properties do not concern us here. We point out, however, that the functors induced by K-schemes are stable and homotopy-invariant by construction.

Roughly speaking, a K-scheme is a selection of elements in each good Fr\'echet algebra in such a way that morphisms and amplifications preserve the selection. An open set $\Omega\subseteq\C$ containing the origin selects, in every good Fr\'echet algebra, those elements whose spectrum is contained in $\Omega$. We can thus associate to each $\Omega$ a K-functor $K_\Omega$; these K-functors are called spectral K-functors. For suitable choices of $\Omega$, one recovers the $K_0$ and $K_1$ functors. We investigate how $K_\Omega$ depends on $\Omega$ and we show, for instance, that conformally equivalent domains yield naturally equivalent spectral K-functors. Finally, we prove the Relative Density Theorem for spectral K-functors.

\subsection{K-schemes and induced K-functors}
A \emph{K-scheme} $S$ associates to each good Fr\'echet algebra $A$ a subset $A_S\subseteq A$ such that the following axioms are satisfied:

$(S_1)$ $0\in A_S$

$(S_2)$ if $a\in \M_p(A)_S$ and $b\in\M_q(A)_S$ then $\diag(a, b) \in \M_{p+q}(A)_S$

$(S_3)$ if $\phi:A\To B$ is a morphism then $\phi(A_S)\subseteq B_S$

\noin Examples will appear shortly. Observe, at this point, that the set of $S$-elements in a good Fr\'echet algebra is invariant under conjugation.

A K-scheme $S$ gives rise to a K-functor $K_S$. We first construct a functor $V_S$ with values in abelian monoids, then we obtain a functor with values in abelian groups via the Grothendieck functor. 

Let $A$ be a good Fr\'echet algebra. The embeddings $\M_n(A)\into \M_{n+1}(A)$, given by $a\mapsto \diag(a, 0)$, restrict to embeddings $\M_n(A)_S\into \M_{n+1}(A)_S$. Put
\[V_S(A)=\bigg(\bigsqcup_{n\geq 1} \M_n(A)_S\bigg)\big/ \sim\]
where the equivalence relation $\sim$ is that of eventual homotopy, that is, $a\in \M_p(A)_S$ and $b\in\M_q(A)_S$ are equivalent if $\diag (a, 0_{n-p})$ and $\diag(b,0_{n-q} )$ are path-homotopic in $\M_n(A)_S$ for some $n\geq p,q$. In other words, $V_S(A)=\varinjlim\; \pi_0\big(\M_n(A)_S\big)$ as a set. The conjugacy-invariance of $S$-elements gives that $V_S(A)$ is an abelian monoid under $[a]+[b]:=[\diag(a, b)]$.

Let $\phi :A\To B$ be a morphism. We then have a monoid morphism $V_S(\phi):V_S(A)\To V_S(B)$ given by $V_S(\phi)([a])=[\phi(a)]$. Here we use $\phi$ to denote each of the amplified morphisms $\M_n(A)\To\M_n(B)$.

So far, we have that $V_S$ is a monoid-valued functor on good Fr\'echet algebras. Let $K_S$ be obtained by applying the Grothendieck functor to $V_S$.  We conclude:

\begin{prop} $K_S$ is a K-functor.
\end{prop}

The K-functor associated to a K-scheme is modeled after the $K_0$ functor, which arises in this way from the \emph{idempotent K-scheme} $A\mapsto\Idem(A)$. One could introduce a ``multiplicative'' version of K-scheme, where the axiom $0\in A_S$ is replaced by $1\in A_S$, and suitably define a corresponding K-functor so that one recovers the $K_1$ functor from the \emph{invertible K-scheme} $A\mapsto A^\times$. However, the difference between the ``additive'' axiom $0\in A_S$ and the ``multiplicative'' axiom $1\in A_S$ is a unit shift, which makes such a ``multiplicative'' K-scheme essentially redundant. Therefore, up to the natural equivalence induced by the unit shift, we may think of the $K_1$ functor as arising form the \emph{shifted invertible K-scheme} $A\mapsto A^\times-1$.

A \emph{morphism of K-schemes} $f:S\To S'$ associates to each good Fr\'echet algebra $A$ a continuous map $f:A_S\To A_{S'}$ such that the following axioms are satisfied:

$(MS_1)$ $f(0)=0$

$(MS_2)$ if $a\in \M_p(A)_S$ and $b\in\M_q(A)_S$ then $\diag(f(a), f(b))=f(\diag(a,b))$ in $\M_{p+q}(A)_{S'}$

$(MS_3)$ if $\phi:A\To B$ is a morphism then the following diagram commutes:
\begin{equation*}
\begin{CD}
A_S @>\phi>> B_S \\
@V f VV @V f VV \\
A_{S'} @>\phi>>B_{S'}
\end{CD}
\end{equation*}

\begin{prop} A morphism of K-schemes $f:S\To S'$ induces a natural transformation of K-functors $K_f:K_S\To K_{S'}$.
\end{prop}

\begin{proof} If suffices to show that $f:S\To S'$ induces a natural transformation $V_f:V_S\To V_{S'}$. First, we show that for any good Fr\'echet algebra $A$ there is a monoid morphism $f_A:V_S(A)\To V_{S'}(A)$. Second, we show that for any morphism $\phi:A\To B$ the following diagram commutes:
\begin{equation*}
\begin{CD}
V_{S} (A)@>V_S(\phi)>> V_S (B) \\
@V f_AVV @V f_B VV \\
V_{S'} (A) @>V_{S'}(\phi) >>V_{S'} (B)
\end{CD}
\end{equation*}
Let $A$ be a good Fr\'echet algebra. Define $f_A:V_S(A)\To V_{S'}(A)$ by $f_A([a])=[f(a)]$. Note that $f_A$ is well-defined: if $a\in M_p(A)_S$ and $b\in\M_q(A)_S$ are eventually homotopic, i.e., $\diag (a, 0_{n-p})$ and $\diag(b,0_{n-q} )$ are path-homotopic in some $\M_n(A)_S$, then 
\[f(\diag (a, 0_{n-p}))=\diag(f(a), 0_{n-p})\]
and 
\[f(\diag(b,0_{n-q} ))=\diag(f(b), 0_{n-q})\] 
are path-homotopic in $\M_n(A)_{S'}$, i.e., $f(a)\in M_p(A)_{S'}$ and $f(b)\in\M_q(A)_{S'}$ are eventually homotopic. Clearly $f_A$ is a monoid morphism: we have $f_A([0])=[f(0)]=[0]$ by $(MS_1)$, while by $(MS_2)$ we have
\begin{eqnarray*}f_A([a]+[b])=f_A[\diag(a,b)]&=&[f(\diag(a,b))]\\&=&[\diag(f(a),f(b))]=[f(a)]+[f(b)]=f_A([a])+f_A([b]).
\end{eqnarray*}
The commutativity of the diagram follows from $(MS_3)$. \end{proof}

Note that, if $f:S\To S'$ and $g:S'\To S''$ are morphisms of K-schemes then $gf: S\To S''$ is a morphism of K-schemes and $K_{gf}=K_g K_f$.

The morphisms of K-schemes $f,g:S\To S'$ are \emph{homotopic} if $f,g: A_S\To A_{S'}$ are homotopic for any good Fr\'echet algebra $A$. If that is the case, then $f$ and $g$ induce the same natural transformation $K_f=K_g: K_S\To K_{S'}$. The K-schemes $S$, $S'$ are \emph{homotopy equivalent} if there are morphisms of K-schemes $f:S\to S'$ and $g:S'\To S$ such that $gf$ is homotopic to $\id_S$ and $fg$ is homotopic to $\id_{S'}$. A K-scheme $S$ is \emph{contractible} if $S$ is homotopy equivalent to the zero K-scheme $A\mapsto \{0_A\}$.

\begin{prop} If the K-schemes $S$ and $S'$ are homotopy equivalent, then the K-functors $K_S$ and $K_{S'}$ are naturally equivalent. In particular, if the K-scheme $S$ is contractible then $K_S$ is the zero functor.
\end{prop}

\subsection{Spectral K-functors: defining $K_\Omega$}
Let $0\in\Omega\subseteq \C$ be an open set. The \emph{spectral K-scheme} $S_\Omega$ associates to each good Fr\'echet algebra $A$ the subset $A_\Omega=\{a: \spec(a)\subseteq \Omega\}$. The axioms are easy to check; for instance, $(S_2)$ follows from
\[\spec_{M_{p+q}(A)}\big(\diag(a,b)\big)=\spec_{M_p(A)}(a)\cup\spec_{M_q(A)}(b)\]
while $(S_3)$ follows from the fact that morphisms are non-increasing on spectra. 

The K-functor associated to the spectral K-scheme $S_\Omega$ is denoted by $K_\Omega$ and is referred to as a spectral K-functor.

\begin{ex} We compute $K_\Omega(\C)$.

Assume, for simplicity, that $\Omega$ has finitely many connected components $\Omega_0, \Omega_1,\dots , \Omega_k$, where $\Omega_0$ is the component of $0$. For $1\leq i\leq k$, let $\#_i(a)$ denote the number of eigenvalues, counted with multiplicity, of $a\in\M_n(\C)$ that lie in $\Omega_i$. Then the eigenvalue-counting map $\#:V_\Omega(\C)\To \N^k$ given by $\#([a])=(\#_1(a),\dots,\#_k(a))$ is well-defined. Visibly, $\#$ is a surjective morphism of monoids. We show $\#$ is  injective. For $1\leq i\leq k$, pick a basepoint $\lambda_i\in\Omega_i$ and think of
\[\lambda(a)=\diag\big(\underbrace{\lambda_1,\dots,\lambda_1}_{\#_1(a)},\dots,\underbrace{\lambda_k,\dots,\lambda_k}_{\#_k(a)}\big)\]
as a normal form for $a\in\M_n(\C)$. It suffices to show that $a$ is eventually homotopic to $\lambda(a)$.  Similar matrices are eventually homotopic; this is true for any K-scheme in fact. Thus $a$ is eventually homotopic to its Jordan normal form. Up to a spectrum-preserving homotopy, one can assume that the Jordan normal form is in fact diagonal. Finally, a homotopy sends all eigenvalues in each $\Omega_i$ to the chosen basepoint $\lambda_i$, and all eigenvalues in $\Omega_0$ to $0$. That is, we can reach a diagonal matrix $\diag(\lambda(a), 0,\dots, 0)$. 

Therefore $V_\Omega(\C)$ is isomorphic to $\N^k$ and, consequently, $K_\Omega(\C)$ is isomorphic to $\Z^k$. \end{ex}

\subsection{Spectral K-functors: dependence on $\Omega$} An open set $\Omega\subseteq \C$ containing $0$ is thought of as a \emph{based} open set. Correspondingly, a holomorphic map $f:\Omega\To \Omega'$ sending $0$ to $0$ is called a \emph{based} holomorphic map.

\begin{prop}\label{holo}
A based holomorphic map $f:\Omega\To \Omega'$ induces a morphism of K-schemes $f_*:S_\Omega\To S_{\Omega'}$, hence a natural transformation of K-functors $K_{f_*}:K_\Omega\To K_{\Omega'}$. In particular, a based conformal equivalence between $\Omega$ and $\Omega'$ induces a natural equivalence between $K_\Omega$ and $K_{\Omega'}$.\end{prop}

\begin{proof}
Let $A$ be a good Fr\'echet algebra. We have  a map $f_*: A_\Omega\To A_{\Omega'}$ given by $f_*(a):=f(a)=\hol_a(f)$, where $\hol_a:\hol(\Omega)\To A$ is the holomorphic calculus of $a\in A_\Omega$. We claim that $f_*$ is continuous. Indeed, let $a_n\To a$ in $A_\Omega$. Pick a (topologically tame) cycle $\gamma$ containing $\spec(a)$ in its interior. Since the set of elements whose spectrum is contained in the interior of $\gamma$ is open, we may assume without loss of generality that all the $a_n$'s have their spectrum contained in the interior of $\gamma$. Then 
\[f_*(a_n)-f_*(a)=\frac{1}{2\pi i}\oint_\gamma f(\lambda)\big((\lambda-a_n)^{-1}-(\lambda-a)^{-1}\big)d\lambda\]
and so $f_*(a_n)\To f_*(a)$ since the integrand converges to $0$. We show $f_*$ is a morphism of K-schemes. Axiom $(MS_1)$ is obvious. For $(MS_2)$ we use the uniqueness of holomorphic calculus. If $\hol_a: \hol(\Omega)\To\M_p(A)$ is the holomorphic calculus for $a$ and $\hol_b:\hol(\Omega)\To\M_q(A)$ is the holomorphic calculus for $b$, then $\diag(\hol_a, \hol_b):\hol(\Omega)\To\M_{p+q}(A)$ given by $g\mapsto \diag(\hol_a(g), \hol_b(g))$
is a holomorphic calculus for $\diag(a, b)$; thus $\hol_{\diag(a, b)} =\diag(\hol_a, \hol_b)$, in particular $f_*(\diag(a,b))=\diag(f_*(a),f_*(b))$. For $(MS_3)$, we use again the uniqueness of holomorphic calculus. If $\hol_a: \hol(\Omega)\To A$ is the holomorphic calculus of $a\in A_\Omega$ then $\phi\circ\hol_a:\hol(\Omega)\To B$ is a holomorphic calculus for $\phi(a)$ hence $\phi\circ\hol_a=\hol_{\phi(a)}$, in particular $\phi(f_*(a))=f_*(\phi(a))$ for all $a\in A_\Omega$.

For the second part, if suffices to check that $(gf)_*=g_* f_*$ for based holomorphic maps $f:\Omega\To\Omega'$ and $g:\Omega'\To\Omega''$. That is, we need $(gf)(a)=g(f(a))$ for all $a\in A_\Omega$. This follows once again by the uniqueness of holomorphic calculus: as both $h\mapsto (hf)(a)$ and $h\mapsto h\big(f(a)\big)$ are holomorphic calculi $\hol(\Omega')\To A$ for $f(a)$, we get $(hf)(a)=h\big(f(a)\big)$ for all $h\in \hol(\Omega')$, in particular for $g$.
\end{proof}

\begin{cor} If $\Omega$ is connected and simply connected then $K_\Omega$ is the zero functor.
\end{cor}

\begin{proof} By the conformal invariance of $K_\Omega$ and the Riemann Mapping Theorem, we may assume that $\Omega$ is either the entire complex plane, or the open unit disk. To show that $K_\Omega$ is the zero functor, it suffices to get $\pi_0(A_\Omega)=0$ for every good Fr\'echet algebra $A$. Indeed, each $a\in A_\Omega$ can be connected to $0_A$ by the path $t\mapsto ta$, path which lies in $A_\Omega$ since $t\Omega\subseteq \Omega$ for $0\leq t\leq 1$. 

In a more sophisticated formulation, $S_\Omega$ is contractible.
\end{proof}

Proposition~\ref{holo} shows that the functor $K_\Omega$ depends only on the conformal type of $\Omega$. But more is true, in fact: $K_\Omega$ depends only on the holomorphic homotopy type of $\Omega$. Roughly speaking, the notion of holomorphic homotopy type is obtained from the usual notion of homotopy type by requiring the maps to be holomorphic rather than continuous.

The based holomorphic maps $f,g:\Omega\To\Omega'$ are \emph{holomorphically homotopic} if there is a family of based holomorphic maps $\{h_t\}_{0\leq t\leq 1}:\Omega\To\Omega'$ such that $h_0=f$, $h_1=g$ and $t\mapsto h_t\in \hol(\Omega)$ is continuous. If that is the case, then the induced morphisms of K-schemes $f_*,g_*:S_\Omega\To S_{\Omega'}$ are homotopic. If there are based holomorphic maps $f:\Omega\To\Omega'$ and $g:\Omega'\To\Omega$ such that $fg$ is holomorphically homotopic to $\id_{\Omega'}$ and $gf$ is holomorphically homotopic to $\id_{\Omega}$, then $\Omega$ and $\Omega'$ are said to be \emph{holomorphic-homotopy equivalent}. If that is the case, then the K-schemes $S_\Omega$ and $S_{\Omega'}$ are homotopy equivalent in the sense described in the previous subsection. Therefore:

\begin{prop} If $\Omega$ and $\Omega'$ are holomorphic-homotopy equivalent, then $K_\Omega$ and $K_{\Omega'}$ are naturally equivalent.
\end{prop}

\subsection{Spectral K-functors: recovering $K_0$ and $K_1$} The $K_1$ functor can be obtained from the shifted invertible K-scheme $A\mapsto A^\times-1$. That is, $K_1$ is naturally equivalent to $K_\Omega$ for $\Omega=\C\setminus\{-1\}$. The $K_0$ functor can be obtained from the idempotent K-scheme $A\mapsto\Idem(A)$. However, the idempotent scheme is not a spectral K-scheme. We now realize the $K_0$ functor as a spectral K-functor.

\begin{prop}\label{K0} $K_0=K_\Omega$ for $\Omega=\C\setminus\{\re=\frac{1}{2}\}$.
\end{prop} 

\begin{proof} Let $\chi$ denote the function defined on $\Omega$ as $\chi=0$ on $\{\re< \frac{1}{2}\}$ and $\chi=1$ on $\{\re> \frac{1}{2}\}$. Consider the holomorphic functions $\{h_t\}_{0\leq t\leq 1}:\Omega\To\Omega$ defined as $h_t=(1-t)\id + t\chi$. We claim that $\{(h_t)_*\}_{0\leq t\leq 1}$ is a strong deformation of $S_\Omega$ to the idempotent K-scheme. It will then follow that $K_\Omega=K_0$.

Let $A$ be a good Fr\'echet algebra. We need to show that $\{(h_t)_*\}_{0\leq t\leq 1}$ is a strong deformation from $A_\Omega$ to the idempotents of $A$. Each $(h_t)_*:A_\Omega\To A_\Omega$ is continuous. The map $t\mapsto h_t$ is $\hol(\Omega)$-continuous, so $t\mapsto (h_t)_*$ is $A$-continuous. At $t=0$, $(h_0)_*=\id_{A_\Omega}$. At $t=1$, $(h_1)_*=\chi_*$ takes idempotent values since $\chi^2=\chi$. Finally, each $(h_t)_*$ acts identically on idempotents. Indeed, it suffices to show that $\chi_*$ acts identically on idempotents. Letting $e$ be an idempotent, we have
\[\chi(e)=\frac{1}{2\pi i}\oint_{\gamma_1} (\lambda-e)^{-1}d\lambda=\frac{1}{2\pi i}\oint_{\gamma_1} \Big(\frac{1-e}{\lambda}+\frac{e}{\lambda-1}\Big)d\lambda=e\]
where $\gamma_1$ is a curve around $1$. \end{proof}

That $K_0$ can be described in terms of elements with spectrum contained in $\C\setminus\{\re=\frac{1}{2}\}$ is discussed in \cite[pp. 193-196]{RLL00}; the above proof is essentially the one given there. It is this alternate picture for $K_0$ that inspired us in defining the $K_\Omega$ groups. Note that $K_\Omega$ is naturally equivalent to $K_0$ whenever $\Omega$ is the disjoint union of two connected and simply connected open subsets of $\C$.

\subsection{Spectral K-functors: the Density Theorem} Finally, we prove the Relative Density Theorem for spectral K-functors. In particular, we obtain the Relative Density Theorem for the usual K-theory. In the case of $K_0$, this proof is more elementary than the proof of Corollary~\ref{quick}, which used the Bott periodicity.
\begin{prop}\label{Kisom}
If $\phi:A\To B$ is a dense and completely relatively spectral morphism between good Fr\'echet algebras, then $\phi$ induces an isomorphism $K_\Omega(A)\simeq K_\Omega(B)$.
\end{prop}

\begin{proof}
Let $\phi:A\To B$ be a dense morphism that is spectral relative to a dense subalgebra $X\subseteq A$. Obviously $A_\Omega\cap X=\phi^{-1}(B_\Omega)\cap X$, so $\phi$ induces a bijection $\pi_0(A_\Omega)\To\pi_0(B_\Omega)$ by Lemma~\ref{key}. Thus, if $\phi:A\To B$ is a dense and completely relatively spectral morphism then $\phi$ induces a bijection $\pi_0(\M_n(A)_\Omega)\To\pi_0(\M_n(B)_\Omega)$ for each $n\geq 1$. It follows that $\phi$ induces a bijection $\varinjlim\;\pi_0(\M_n(A)_\Omega)\To\varinjlim\;\pi_0(\M_n(B)_\Omega)$. In other words, $V_\Omega(\phi):V_\Omega(A)\To V_\Omega(B)$ is a bijection. We conclude that $\phi$ induces an isomorphism $K_\Omega(A)\To K_\Omega(B)$.
\end{proof}







\section{Applications}\label{App}

\subsection{Spectral equivalence}
Let  $A$ and $B$ be Banach algebra completions of an algebra $X$. Say that $A$ and $B$ are \emph{spectrally equivalent over $X$} if $\spec_A(x)=\spec_B(x)$ for all $x\in X$, and \emph{completely spectrally equivalent over $X$} if $\M_n(A)$ and $\M_n(B)$ are spectrally equivalent over $\M_n(X)$ for each $n\geq 1$.

\begin{prop} The higher connected stable ranks and the property of being finite are invariant under spectral equivalence. Furthermore, K-theory and the property of being stably finite are invariant under complete spectral equivalence.
\end{prop}

\begin{proof} Let  $A$, $B$ be Banach algebra completions of $X$. Let $\overline{X}$ be the Banach algebra obtained by completing $X$ under the norm $\|x\|:=\|x\|_A+\|x\|_B$. Then $r_{\overline{X}}(x)=\max\{r_A(x),r_B(x)\}$ for all $x\in X$. Now, if $A$ and $B$ are spectrally equivalent over $X$ then $r_{\overline{X}}(x)=r_A(x)=r_B(x)$ for all $x\in X$. We obtain dense and relatively spectral morphisms $\overline{X}\To A$, $\overline{X}\To B$ by extending the inclusions $X\into A$, $X\into B$. Thus $\csr_k(A)=\csr_k(\overline{X})=\csr_k(B)$, and $A$ is finite if and only if $\overline{X}$ is finite if and only if $B$ is finite.

If $A$ and $B$ are completely spectrally equivalent over $X$ then we obtain morphisms $\overline{X}\To A$, $\overline{X}\To B$ that are completely spectral relative to $X$. Hence $K_*(A)\simeq K_*(\overline{X})\simeq K_*(B)$, and $A$ is stably finite if and only if $\overline{X}$ is stably finite if and only if $B$ is stably finite.
\end{proof}

Two algebras that are spectrally equivalent need not be connected by a morphism. Consider, however, the following $*$-context: $A$ is a Banach $*$-algebra, $B$ is a $\Cstar$-algebra, and $X$ is a dense $*$-subalgebra of $A$ and $B$. The argument used in the proof of Lemma~\ref{isom} shows the following: if $A$ and $B$ are spectrally equivalent over $X$, then there is a dense and relatively spectral morphism $A\To B$ extending the identity on $X$.

A criterion for complete spectral equivalence is given in the next section.

\subsection{Subexponential control}\label{subexp}
A \emph{weight} on a group $\G$ is a map $S\mapsto \omega(S)\in [1,\infty)$ on the finite nonempty subsets of $\G$ which is non-decreasing, i.e., $S\subseteq S'$ implies $\omega(S)\leq\omega(S')$. A weight $\omega$ is \emph{subexponential} if $\omega(S^n)^{1/n}\To 1$ for all $S$.

\begin{prop}\label{control}
Let $A\G$ and $B\G$ be Banach algebra completions of a group algebra $\C\G$. Assume that there are constants $C, C'>0$ and subexponential weights $\omega$, $\omega'$ such that
\[\|a\|_B\leq C\omega(\supp\; a)\|a\|_A, \quad \|a\|_A\leq C'\omega'(\supp\; a)\|a\|_B\]
for all $a\in \C\G$. Then $A\G$ and $B\G$ are completely spectrally equivalent over $\C\G$.
\end{prop}

\begin{proof} Let $k\geq 1$. For all $(a_{ij})\in \M_k(\C\G)$ we have analogous subexponential estimates in the matrix algebras $\M_k(A\G)$ and $\M_k(B\G)$:
\[\|(a_{ij})\|_B\leq C\omega\big(\supp\; (a_{ij})\big)\|(a_{ij})\|_A, \quad \|(a_{ij})\|_A\leq C'\omega'\big(\supp\; (a_{ij})\big)\|(a_{ij})\|_B\]
As $\omega\big(\supp\; (a_{ij})^n\big)\leq \omega\big((\supp\; (a_{ij}))^n\big)$, we get $r_A\big((a_{ij})\big)=r_B\big((a_{ij})\big)$ for all $(a_{ij})\in \M_k(\C\G)$. \end{proof}

\begin{ex}\label{subgrowth} Let $\G$ be a group of subexponential growth. We claim that $\ell^1\G\into\Cred\G$ is completely spectral relative to $\C\G$. Indeed, let  $\omega$ be the subexponential weight on $\G$ given by $\omega(S)=\sqrt{\mathrm{vol}\; B(S)}$, where $B(S)$ denotes the ball centered at the identity that circumscribes $S$. For all $a\in \C\G$ we have $\|a\|\leq \|a\|_1$ and $\|a\|_1\leq\omega(\supp\; a)\|a\|_2\leq\omega(\supp\; a)\|a\|$, i.e., we are in the conditions of Proposition~\ref{control}. We infer that $K_*(\ell^1\G)\simeq K_*(\Cred\G)$.
\end{ex}

Consider, at this point, the following
\begin{conj}
For any discrete countable group $\G$, the inclusion $\ell^1\G\into\Cred\G$ induces an isomorphism $K_*(\ell^1\G)\simeq K_*(\Cred\G)$.
\end{conj}
Let us refer to this statement as the BBC conjecture, for it connects the Bost conjecture with the Baum-Connes conjecture. We find the BBC conjecture to be a natural question on its own, outside of the Bost-Baum-Connes context. One knows, by combining results of Higson and Kasparov on the Baum-Connes side with results of Lafforgue on the Bost side, that the BBC conjecture holds for groups with the Haagerup property. On the other hand, the spectral approach allowed us to verify the BBC conjecture for groups of subexponential growth. The severe limitations of the spectral approach are made evident by this comparison.

\begin{ex} Let $\G$ be a finitely generated group satisfying the Rapid Decay property, i.e., there are constants $C,d>0$ such that $\|a\|\leq C\|a\|_{2,d}$ for all $a\in \C\G$. The weighted $\ell^2$-norm $\|\cdot\|_{2,d}$ is given by $\big\|\sum a_g g\big\|_{2,d}=\sqrt{\sum |a_g|^2 (1+|g|)^{2d}}$, where $| \cdot |$ denotes the word-length. Examples of groups satisfying the Rapid Decay property include free groups \cite{Haa79} and, more generally, hyperbolic groups \cite{dHa88}, groups of polynomial growth \cite{Jol90}, and many other groups \cite{CR05}, \cite{DS05}, \cite{Laf00}, \cite{RRS}.

Consider the weighted $\ell^2$-space $\ell^2_s\G=\big\{\sum a_g g:\; \big\|\sum a_g g\big\|_{2,s}<\infty\big\}$. For $s>d$, $\ell^2_s\G$ is a Banach subalgebra of $\Cred\G$ (\cite[Prop.1.2]{Laf00}). Lafforgue's goal is the isomorphism $K_*(\ell^2_s\G)\simeq K_*(\Cred\G)$, so he shows that $\ell^2_s\G$ is a spectral subalgebra of $\Cred\G$.

Alternatively, one can adopt the relative perspective. Let $\omega_s$ be the subexponential weight on $\G$ given by $\omega_s(S)=\big(1+ R(S)\big)^s$, where $R(S)$ is the radius of the ball centered at the identity that circumscribes $S$. For $s>d$, we have $\|a\|\leq C\|a\|_{2,s}$ and $\|a\|_{2,s}\leq \omega_s(\supp\; a)\|a\|_2\leq \omega_s(\supp\; a)\|a\|$ for all $a\in \C\G$. From Proposition~\ref{control} it follows that $K_*(\ell^2_s\G)\simeq K_*(\Cred\G)$.
\end{ex}

\subsection{$\Sigma_1$-groups}\label{Sigma1}
We consider groups $\G$ for which the inclusion $\ell^1\G\into\Cred\G$ is spectral relative to $\C\G$. Let us refer to such groups $\G$ as \emph{$\Sigma_1$-groups}. That is, $\G$ is a $\Sigma_1$-group if $r_{\ell^1\G}(a)\leq \|a\|$ for all $a\in \C\G$. The $\Sigma_1$ condition is the $\ell^1$ analogue of the ``$\ell^2$-spectral radius property'' discussed in \cite[Def.1.2 ii) \& Sec.3]{DdH99}. The major difference is that the ``$\ell^2$-spectral radius property'' is satisfied by Rapid Decay groups, e.g. by hyperbolic groups, whereas $\Sigma_1$-groups are necessarily amenable:

\begin{prop}\label{Sigma} We have:

a) $\Sigma_1$ is closed under taking subgroups;

b) $\Sigma_1$ is closed under taking directed unions;

c) $\Sigma_1$-groups are amenable;

d) $\Sigma_1$-groups do not contain $FS_2$, the free semigroup on two generators;

e) groups of subexponential growth are $\Sigma_1$-groups.
\end{prop}

\begin{proof}
a) An embedding $\Lambda\into\G$ induces isometric embeddings $\ell^1\Lambda\into\ell^1\G$, $\Cred\Lambda\into \Cred\G$.

b) Let $\G$ be the directed union of the $\Sigma_1$-groups $(\G_i)_{i\in I}$. If $a\in \C\G$, then $a\in\C\G_i$ for some $i$ and we know that $r_{\ell^1\G_i}(a)=r_{\Cred\G_i}(a)$ can also be read as $r_{\ell^1\G}(a)=r_{\Cred\G}(a)$.

c) Let $\G$ be a $\Sigma_1$-group. The inclusion $\ell^1\G\into\Cred\G$ factors as $\ell^1\G\into\Cstar\G\onto\Cred\G$. It follows that $\Cstar\G\onto\Cred\G$ is spectral relative to $\C\G$, so necessarily $\Cstar\G\onto\Cred\G$ is an isomorphism, i.e., $\G$ is amenable.

Here is another argument. By a), it suffices to show that $\G$ is amenable in the case $\G$ is finitely generated, say by a finite symmetric set $S$. For $a\in \R_+\G$ we have $\|a^n\|_1=\|a\|^n_1$, so $r_{\ell^1\G}(a)=\|a\|_1$. We obtain $r_{\Cred\G}(\chi_S)=r_{\ell^1\G}(\chi_S)=\|\chi_S\|_1=\#S$. By Kesten's criterion, $\G$ is amenable.

d) Let $\G$ be a $\Sigma_1$-group and assume, on the contrary, that $x,y\in\G$ generate a free subsemigroup. If the support of $a\in \C\G$ generates a free subsemigroup then $\|a^n\|_1=\|a\|^n_1$, hence $r_{\ell^1\G}(a)=\|a\|_1$. Together with the $\Sigma_1$ condition, this gives $\|a\|=\|a\|_1$ for those $a\in \C\G$ whose support generates a free subsemigroup. Consider now $a=1+ix+ix^{-1}\in\C\G$. On one hand, the support of $(yx)a$ generates a free subsemigroup, so $\|a\|=\|(yx)a\|=\|(yx)a\|_1=\|a\|_1=3$. On the other hand, we have $\|a\|^2=\|aa^*\|=\|3+x^2+x^{-2}\|\leq 5$, a contradiction.

e) This was proved in Example~\ref{subgrowth}.
\end{proof}

Thus $\Sigma_1$-groups include the subexponential groups and are included among amenable groups without free subsemigroups. Note that there are examples, first constructed by Grigorchuk, of amenable groups of exponential growth that do not contain free subsemigroups.

\begin{rem} If one considers the groups $\G$ for which the inclusion $\ell^1\G\into\Cstar\G$ is spectral relative to $\C\G$, then the analogues of parts a), b), d), e) of Proposition~\ref{Sigma} still hold.
\end{rem}

After Lafforgue, we say that a Banach algebra completion $A\G$ of $\C\G$ is an \emph{unconditional completion} if the norm $\|\cdot\|_A$ of $A\G$ satisfies $\||a|\|_A=\|a\|_A$ for all $a\in\C\G$, where $|a|$ denotes the pointwise absolute value. The simplest example of unconditional completion is $\ell^1\G$. For groups with the Rapid Decay property, $\ell^2_s\G$ ($s\gg 1$) is an unconditional completion. In light of \cite{Laf02}, one is interested in finding unconditional completions $A\G$ having the same K-theory as $\Cred\G$. One way to achieve $K_*(A\G)\simeq K_*(\Cred\G)$ would be to have $A\G$ and $\Cred\G$ completely spectrally equivalent over $\C\G$. For amenable groups, however, one cannot do better than $\Sigma_1$-groups:

\begin{prop} Let $\G$ be amenable. If there is an unconditional completion $A\G$ with the property that $A\G$ and $\Cred\G$ are spectrally equivalent over $\C\G$, then in fact $\ell^1\G$ and $\Cred\G$ are spectral over $\C\G$, i.e., $\G$ is a $\Sigma_1$-group.
\end{prop}

\begin{proof} We first show that $r_{\ell^1\G}(a)=r_{\Cred\G}(a)$ whenever $a\in\C\G$ is self-adjoint. The spectral equivalence of $A\G$ and $\Cred\G$ yields $\|a\|=r_{\Cred\G}(a)=r_{A\G}(a)\leq \|a\|_A$ for all self-adjoint $a\in\C\G$. Since $\G$ is amenable, we get $\|a\|_1=\||a|\|\leq \||a|\|_A=\|a\|_A$ for all self-adjoint $a\in\C\G$, so $\|a\|\leq \|a\|_1\leq \|a\|_A$ for all self-adjoint $a\in\C\G$. It follows that $r_{\ell^1\G}(a)=r_{\Cred\G}(a)$ for all self-adjoint $a\in\C\G$.

Next, one needs to adapt the proof of Proposition~\ref{radius} to the following $*$-version:

\begin{quote}
Let $\phi:A \To B$ be a dense $*$-morphism between Banach $*$-algebras, and let $X\subseteq A$ be a dense $*$-subalgebra. If $r_B(\phi(x))=r_A (x)$ for all self-adjoint $x\in X$, then $\phi:A \To B$ is spectral relative to $X$.
\end{quote}

Indeed, let $x\in X$ with $\phi(x)$ invertible in $B$; we show $x$ invertible in $A$. Pick $(x_n)\subseteq X$ such that $\phi(x_n)\To\phi(x)^{-1}$. Then $\phi((xx_n)(xx_n)^*)\To 1$, so 
\[r_A(1-(xx_n)(xx_n)^*)=r_B\big(\phi(1-(xx_n)(xx_n)^*) \big)\To 0.\]
 It follows that $(xx_n)(xx_n)^*$ is invertible in $A$ for large $n$. Similarly, $(x_nx)^*(x_nx)$ is invertible in $A$ for large $n$. We conclude that $x$ is invertible in $A$. \end{proof}

The above proof becomes even simpler if one makes the (natural) assumption that the unconditional completion $A\G$ is a Banach $*$-algebra. One then argues, as in the proof of Lemma~\ref{isom}, that $\|a\|\leq \|a\|_A$ for all $a\in\C\G$.


\end{document}